\documentclass[12pt,leqno,french]{amsart}

\def\opn#1#2{\def#1{\mathop{\kern0pt\fam0#2}\nolimits}}

\opn\arg{arg} \opn\DD{\mathcal D} \opn\GCD{GCD}
\opn\Id{Id} \opn\im{Im} \opn\rad{rad}
\opn\sgdeg{sgdeg}  \opn\Frac{Frac} \opn\Spec{Spec}

\usepackage[latin1]{inputenc}
\usepackage{babel}
\usepackage{amsmath}
\usepackage{amsfonts}
\usepackage{amscd}
\usepackage{amssymb}

\begin{document}

\begin{center}{\large \bf Sur le Th\'eor\`eme Principal de Zariski en G\'eom\'etrie Alg\'ebrique et G\'eom\'etrie Analytique}\\[0.2cm]
{Kossivi Adjamagbo \\ Universit\'e Paris 6 - Case 172 - Institut de Mathématiques de Jussieu\\ 4, place Jussieu, 75252 PARIS CEDEX 05\\ adja@math.jussieu.fr}
\end{center}

\vspace*{.2cm}

{\tiny {\bf Résumé :} Nous comblons une  lacune étonnante de la Géométrie Analytique Complexe en prouvant l'analogue du Théorème Principal de Zariski dans cette géométrie, c'est-à-dire en prouvant que toute application holomorphe d'un espace analytique irreductible dans un espace analytique normal et irreductible est un plongement ouvert si et seulement si toutes ses fibres sont discrètes et si elle induit une application biméromorphe sur son image. Nous prouvons plus généralement le ``Théorème Principal de Zariski Généralisé pour les espaces analytiques'', qui affirme qu'une application holomorphe d'un espace analytique irreductible dans un espace analytique irreductible et localement irreductible est un plongement ouvert si et seulement si elle est plate et induit une application biméromorphe sur son image. Gr\^ace au ``critère analytique de régularité'' de Serre-Samuel dans GAGA \cite{serre} et au ``Principe de Lefschetz'', nous en d\'eduisons enfin le ``Théorème Principal de Zariski Généralisé pour les variétés algébriques de caractéristique nulle'', qui affirme qu'un morphisme d'une telle variété irreductible dans une autre unibranche est une immersion ouverte si et seulement s'il est birationnel et plat.}

\vspace*{.2cm}

{\tiny {\bf Abstract : On Zariski Main Theorem in Algebraic Geometry and Analytic Geometry.} We fill a surprising gap of Complex Analytic Geometry by proving the analogue of Zariski Main Theorem in this geometry, i.e. proving that an holomorphic map from an irreducible analytic space to a normal irreducible one is an open embedding if and only if all its fibers are discrete and it induces a bimeromorphic map on its image. We prove more generally the ``Generalized Zariski Main Theorem for analytic spaces'', which claims that an holomorphic map from an irreducible analytic space to a irreducible locally irreducible one is an open embedding if and only if it is flat and induces a bimeromorphic map on its image. Thanks to the ``analytic criterion of regularity'' of Serre-Samuel in GAGA \cite{serre} and to ``Lefschetz Principle'', we finally deduce the ``Generalized Zariski Main Theorem for algebraic varieties of characteristical zero'', which claims that a morphism from such an irreducible variety to an irreducible unibranch one is an open immersion if and only if it is birational and flat.}

\vspace*{.5cm}

\noindent
{\bf D\'efinition.}

1) Une application méromorphe (au sens de Grauert-Remmert, voir par exemple \cite{detgra}, p. 211) entre des espaces complexes est dite ``biméromorphe'' si elle induit une application biholomorphe entre les complémentaires de  sous-ensembles analytiques fermés d'intérieur vide de ces espaces (voir par exemple \cite{peternell}, p. 305).

2) Une application holomorphe entre des espaces complexes est dite ``plate'' si, en tout point l'espace de départ, elle induit un homomorphisme plat entre l'anneau local de l'espace d'arrivée en l'image de ce point par l'application et l'anneau local de l'espace de départ en ce point.

\noindent
{\bf Théorème 1 (Théorème Principal de Zariski pour les espaces analytiques)}. Une application holomorphe d'un espace analytique irreductible dans un espace analytique normal et irreductible de même dimension est un plongement ouvert si et seulement si toutes ses fibres sont discrètes et si elle induit une application biméromorphe entre son espace de départ et son image qui est alors un ouvert de l'espace d'arrivée.

\begin{proof}

Soit $f : X \longrightarrow X'$ une telle application. La condition indiquée étant trivialement nécessaire, prouvons qu'elle est suffisante. Supposons la donc. D'après le Lemme-clé ci-après, $f$ est alors une application ouverte, en particulier l'image de $f$ est un ouvert de $X'$. Or, d'après le Lemme-bis suivant, $f$ est injective. Le même Lemme-clé permet alors de conclure. 

\end{proof}

\noindent
{\bf Lemme-clé (critère d'ouverture et de plongement ouvert pour les applications holomorphes)}. Toute application holomorphe à fibres discrètes entre deux espaces analytiques irreductibles de même dimension dont le second est localement irréductible est une application ouverte. Par suite, si en plus cette application est injective, alors elle induit un homéomorphisme entre son espace de départ et son image. Si en plus l'espace d'arrivée est normal, alors cette application est un plongement ouvert.

\begin{proof}

Soient $f : X \longrightarrow X'$ une telle application, $p : X'' \longrightarrow X$ la projection canonique d'une normalisation de $X$ dans $X$ (voir par exemple \cite{remmert}, Th. 14.9, p. 88) et $g = f \circ p$. Puisque $p$ est à fibres finies, $g$ est donc à fibres discrètes. $X''$ étant irréductible et localement irreductible de même dimension que $X$, il résulte de \cite{remmert}, 10.5, p. 64, que $g$ est ouvert. Il résulte alors de la surjectivité de $p$ que $f$ aussi est ouvert. La conclusion résulte alors de \cite{remmert}, 15.4 et 15.5, p. 93.

\end{proof}

\noindent
{\bf Remarque 1}. Le Lemme-clé est l'analogue fidèle en géométrie analytique du Lemme-clé 14.4.1.3 de EGA \cite{grothendieck3} qui affirme, compte tenu de \cite{grothendieck2}, Prop. 5.4.1, que tout morphisme d'un schéma irréductible dans un schéma irréductible unibranche est ouvert.

\noindent
{\bf Lemme-bis (critère d'injectivité pour les applications holomorphes biméromorphes).} Toute application holomorphe biméromorphe et ouverte est injective.

\begin{proof}

1) Soit $f : X \longrightarrow X'$ une telle application. Notons $F$ (resp. $F'$) un sous-ensemble analytique fermé de $X$ (resp. $X'$) d'intérieur vide tel que $f$ induise une application biholomorphe de $X - F$ dans $X' - F'$.

2) Supposons $f$ non injective et considérons des points distincts $a$ et $b$ de $X$ tels que $c = f(a) = f(b)$. D'après le Lemme-ter suivant, il existe donc des voisinages ouverts disjoints $V_a$ et $V_b$ de $a$ et de $b$ dans $X$  tels que $V_c = f(V_a) = f(V_b)$ soit un voisinage ouvert de $c$ dans $X'$ contenant un sous-ensemble analytique $F''$ dans lequel sont contenus $f(F \cap V_a)$ et $f(F \cap V_b)$. 

3) D'après les choix de $F'$ et $F''$, $V_c - F' - F''$ est non vide. Soit donc $c'$ un de ses points. D'après les choix de $V_a$, $V_b$ et $V_c$, il existe donc $a'$ dans $V_a - F$ et $b'$ dans $V_b - F$ tels que $c' = f(a') = f(b')$. Puisque $V_a$ et $V_b$ sont disjoints, $a'$ et $b'$ sont distincts, ce qui contredit l'injectivité de la restriction de $f$ à $X - F$.

\end{proof}

\noindent
{\bf Lemme-ter (sur les images locales des applications holomorphes).} Pour toute application holomorphe entre des espaces complexes, tout point de son espace de départ admet un voisinage ouvert dont l'image par cette application est contenue dans un sous-ensemble analytique d'un voisinage de l'image de ce point et dont la dimension est au plus égale à celle de l'espace de départ.

\begin{proof}

 Soient $f : X \longrightarrow X'$ une telle application , $a$ un point de $X$, $V$ un voisinage ouvert de $a$ que l'on peut supposer plongée dans un espace projectif complexe $P$, $W$ un voisinage ouvert de $f(a)$ que l'on peut supposer plongé également dans $P$, $f'$ l'application de $V$ dans $W$ induite par $f$, $\Gamma$ le graphe de $f'$, $\bar{\Gamma}$ la fermeture de $\Gamma$ dans $P^2$, et $A$ la seconde projection de $\bar{\Gamma}$ sur $P$. D'après le Théorème de Remmert sur l'image d'une application holomorphe propre (voir par exemple \cite{lowa}, p. 290), $A$ est un sous-ensemble analytique de $P$ de dimension au plus celle de $X$. $A\cap W$ est donc un sous-ensemble analytique de $W$ contenant l'image de $V$ par $f$ et de dimension au plus celle de $X$.
 
\end{proof} 

\noindent
{\bf Théorème 2 (Théorème Principal de Zariski Généralisé pour les espaces analytiques).} Une application holomorphe d'un espace analytique irreductible dans un espace analytique irreductible et localement irreductible est un plongement ouvert si et seulement si elle est plate et induit un application biméromorphe entre son espace de départ et son image qui est alors un ouvert de l'espace d'arrivée.

\begin{proof}

1) Soit $f : X \longrightarrow X'$ une telle application.

2) D'après l'hypothèse de biméromorphie, $X$ et $X'$ ont la même dimension. Il résulte alors de la formule bien connue sur les dimension des fibres d'un morphisme plat, et plus généralement ouvert, d'anneaux locaux noethériens que toutes les fibres de $f$ sont discrètes (voir par exemple \cite{grothendieck1} et Th. 2.4.6, \cite{grothendieck3}, Th. 14.2.1 ou \cite{matsumura}, Th. 15.1, p. 116).

3) Par suite, d'après le Lemme-clé, $f$ est une application ouverte, en particulier l'image de $f$ est un ouvert de $X'$.

4) D'autre part, d'après le Lemme-bis, $f$ est injective. Toujours d'après le Lemme-clé, $f$ induit  donc un homéomorphisme entre $X$ et $f(X)$.

5) Considérons maintenant un point $x$ de $X$ et notons $A$ l'anneau local de $X'$ en $f(x)$, $B$ celui de $X$ en $x$. $X$ étant localement irreductible, l'anneau $A$ est donc intègre. Il admet donc un corps de fractions $C$. 

6) D'après la platitude de $f$, nous pouvons supposer que $A$ est un sous-anneau de $B$. $A$ et $B$ étant des anneaux locaux, $B$ est donc fidèlement plat sur $A$ (voir par exemple \cite{altman}, Ch. V, Prop. 1.6, p. 84). Par suite, $B/A$ est un $A$-module plat, donc sans $A$-torsion (voir par exemple \cite{altman}, Ch. V, Th. 1.9, p. 85).  L'application canonique du $A$-module $B/A$ dans son localisé total $M = (B/A) {\otimes}_A C$ est donc injective, puisque $A$ est intègre. Or d'après la biméromorphie de $f$, $M$ est nul car $f$ induit un isomorphisme entre les faisceaux des fonctions méromorphes sur $X$ et $f(X)$ (voir par exemple \cite{peternell}, Cor. 6.8, p. 305). On en déduit que $A= B$.

7) Compte tenu de 4), on conclut donc que $f$ induit un isomorphisme entre les faisceaux structuraux de $X$ et $f(X)$, c'est-à-dire que $f$ est un plongement ouvert comme désiré.

\end{proof}

\noindent
{\bf Théorème 3 (sur le lien entre le Théorème Principal de Zariski et le Théorème Principal de Zariski Généralisé pour les espaces analytiques).} Le Théorème Principal de Zariski pour les espaces analytiques résulte du Théorème Principal de Zariski Généralisé pour les espaces analytiques.

\begin{proof}

1) Supposons donc ce dernier théorème et considérons $f : X \longrightarrow X'$ une application vérifiant les hypothèses de ce premier théorème, $p : X'' \longrightarrow X$ la projection canonique d'une normalisation de $X$ dans $X$, $g = f \circ p$, $f'$ l'application de $X$ dans $f(X)$ induite par $f$ et $g'$ l'application de $X''$ dans $f(X)$ induite par $g$. 

2) D'après le Lemme-clé, $f(X)$ est un ouvert de $X'$. Puisque $f'$ et $p$ sont biméromorphes à fibres discrètes, d'après le Lemme-clé et le Lemme-ter, tout point $x$ de $X''$ admet un voisinage ouvert $V_x$ que l'on peut supposer plongé dans un espace affine complexe $E$ tel que $g'$ induise une application biméromorphe et à fibres discrètes $g'_x$ de $V_x$ sur l'ouvert $g'(V_x)$ de $X'$. Pour tout $x$ de $X''$, notons $g''_x$ l'application de $V_x - g'^{-1}(S(g'(V_x)))$ dans $g'(V_x) - S(g'(V_x))$, où $S(g'(V_x))$ désigne le lieu singulier de la variété normale $g'(V_x)$.  

3) Considérons maintenant un point $x$ de $X''$. La variété normale $V_x$ étant une variété de Cohen-Macaulay et la variété $g'(V_x) - S(g'(V_x))$ étant non singulière de même dimension que $V_x$, l'application à fibres discrètes $g''_x$ est biméromorphe et plate, en vertu du critère de platitude en termes de fibres pour les homomorphismes d'anneaux locaux (voir par exemple \cite{matsumura}, Th. 23.1, p. 179). D'après le Théorème Principal de Zariski Généralisé pour les espaces analytiques, $g''_x$ est donc biholomorphe. Il résulte alors du théorème de prolongement de Riemann pour les espaces complexes normales que $g'_x$ est biholomorphe (voir par exemple \cite{remmert}, Th. 13.6, p. 81). 

4) $g$ est donc plat. Toujours d'après le Théorème Principal de Zariski Généralisé pour les espaces analytiques, $g'$ est donc biholomorphe. On en conclut que $p$ et $f'$ sont bijectives et que $f'^{-1} = p \circ g'^{-1}$ est holomorphe, ce qui prouve que $f$ est un plongement ouvert et établit donc le Théorème Principal de Zariski pour les espaces analytiques.

\end{proof}

\noindent
{\bf Théorème 4 (Th\'eor\`eme Principal de Zariski Généralisé pour les variétés algébriques complexes).} Un morphisme d'une variété algébrique complexe irréductible dans une variété algébrique complexe irréductible et  unibranche est une immersion ouverte si et seulement s'il est birationnel et plat. 

\begin{proof}

L'énoncé résulte du Théorème Principal de Zariski Généralisé pour les espaces analytiques grâce au ``critère analytique de régularité'' de Serre-Samuel dans \cite{serre}, Prop. 9, p. 13-15 (et la note en bas de la page 13), sachant qu'un homomorphisme d'anneaux locaux noethériens est plat si et seulement si l'homomorphisme induit entre les complétés de ces anneaux l'est (voir par exemple \cite{altman}, Ch. V, Prop. 3.3, p. 93).

\end{proof}

\noindent
{\bf Corollaire 1 (Théorème Principal de Zariski pour les variétés algébriques complexes \cite{grothendieck1}, Cor. 4.4.9, p. 137 et \cite{grothendieck3}, Cor. 8.12.10, p. 48, \cite{mumford}, Ch. III.9, p. 288).} Un morphisme d'une variété algébrique complexe irréductible dans une variété algébrique complexe irréductible normale est une immersion ouverte si et seulement s'il est birationnel et à fibres finies.

\begin{proof}

Soit $f : X \longrightarrow X'$ un morphisme birationnel et à fibres finies. D'après la remarque ci-dessus, $f(X)$ est donc un ouvert de Zariski de $X'$. La conclusion résulte alors du Théorème Principal de Zariski pour les espaces analytiques, grâce au ``critère analytique de régularité'' de Serre-Samuel.

\end{proof}

\noindent
{\bf Corollaire 2 (Théorème Principal de Zariski pour les variétés algébriques de caractéristique nulle, idem).} Un morphisme d'une variété algébrique irréductible sur un corps algébriquement clos de caractéristique nulle dans une variété algébrique irréductible normale sur ce corps est une immersion ouverte si et seulement s'il est birationnel et à fibres finies.

\begin{proof}

En effet, un tel énoncé résulte du précédent par ``le principe de Lefschetz'' fondé sur la possibilité de plonger dans le corps des nombres complexes tout corps de nombres, c'est-à-dire tout corps commutatif qui est une extension finie du corps des nombres rationnels.

\end{proof}

\noindent
{\bf Corollaire 3 (Théorème Principal de Zariski Généralisé pour les variétés algébriques de caractéristique nulle).} Un morphisme d'une variété algébrique irréductible sur un corps algébriquement clos de caractéristique nulle dans une variété algébrique irréductible unibranche sur ce corps est une immersion ouverte si et seulement s'il est birationnel et plat.

\begin{proof}

En effet, un tel énoncé résulte du Théoème 4 par ``le principe de Lefschetz''. 

\end{proof}

\noindent
{\bf Théorème 5 (sur le lien entre le Théorème Principal de Zariski et le Théorème Principal de Zariski Généralisé pour les variétés algébriques de caractéristique nulle).} Le Théorème Principal de Zariski pour les variétés algébriques de caractéristique nulle résulte du Théorème Principal de Zariski Généralisé pour de telles variétés.

\begin{proof}

1) Supposons donc ce dernier théorème et considérons $f : X \longrightarrow X'$ un morphisme vérifiant les hypothèses de ce premier théorème, $p : X'' \longrightarrow X$ la projection canonique d'une normalisation de $X$ dans $X$, $g = f \circ p$, $f'$ l'application de $X$ dans $f(X)$ induite par $f$ et $g'$ l'application de $X''$ dans $f(X)$ induite par $g$. 

2) D'après la remarque 1, $f(X)$ est un ouvert de Zariski de $X'$. Puisque $f'$ et $p$ sont birationnelles à fibres finies, toujours d'après cette remarque, tout point $x$ de $X''$ admet un voisinage ouvert affine $V_x$  tel que $g'$ induise un morphisme birationnel et à fibres finies $g'_x$ de $V_x$ sur l'ouvert de Zariski $g'(V_x)$ de $X'$. Pour tout $x$ de $X''$, notons $g''_x$ l'application de $V_x - g'^{-1}(S(g'(V_x)))$ dans $V'_x - S(g'(V_x))$, où $S(g'(V_x))$ désigne le lieu singulier de la variété normale $g'(V_x)$.  

3) Considérons maintenant un point $x$ de $X''$. La variété normale $V_x$ étant une variété de Cohen-Macaulay et la variété $g'(V_x) - S(g'(V_x))$ étant non singulière de même dimension que $V_x$, le morphisme à fibres finies $g''_x$ est birationnel et plat, en vertu du critère de platitude en termes de fibres pour les homomorphismes d'anneaux locaux. D'après le Théorème Principal de Zariski Généralisé pour les variétés algébriques de caractéristique nulle, $g''_x$ est donc un isomorphisme. Il résulte alors du théorème de prolongement des fonctions régulières sur les variétés algébriqes normales que $g'_x$ est un isomorphisme (voir par exemple \cite{iitaka}, Ch. 2, Th. 1.5, p. 124). 

4) $g$ est donc plat. Toujours d'après le Théorème Principal de Zariski Généralisé pour les variétés algébriques de caractéristique nulle, $g'$ est donc un isomorphisme. On en conclut que $p$ et $f'$ sont bijectives et que $f'^{-1} = p \circ g'^{-1}$ est un morphisme, ce qui prouve que $f$ est un plongement ouvert et établit donc le Théorème Principal de Zariski pour les variétés algébriques de caractéristique nulle.

\end{proof}


\vspace*{.5cm}










\vspace{0.1cm}

\begin{thebibliography}{10}

\bibitem{altman}
A. Altman, S. Kleiman,
\newblock  {\em Introduction to Grothendieck Duality Theorem},
\newblock  Lect. Notes in Math., No. {\bf 146}, Springer-Verlag, 1970.

\bibitem{detgra}
G. Dethloff, H. Grauert,
\newblock  {\em Seminormal Complex Spaces},
\newblock  in Several Complex Variables VII, H. Grauert, Th. Peternell, R. Remmert (Eds), Encyclopedia of Mathematical Sciences {\bf 74} (1994), p. 183-220.


\bibitem{grothendieck1}
A. Grothendieck, J. Dieudonné,
\newblock  {\em Eléments de Géométrie Algébrique, III, Première Partie},
\newblock  I.H.E.S., Public. Math., No. {\bf 11}(1961).

\bibitem{grothendieck2}
A. Grothendieck, J. Dieudonné,
\newblock  {\em Eléments de Géométrie Algébrique, IV, Seconde Partie},
\newblock  I.H.E.S., Public. Math., No. {\bf 24}(1965).

\bibitem{grothendieck3}
A. Grothendieck, J. Dieudonné,
\newblock  {\em Eléments de Géométrie Algébrique, IV, Troisième Partie},
\newblock  I.H.E.S., Public. Math., No. {\bf 28}(1966).

\bibitem{iitaka}
Shigeru Iitaka,
\newblock  {\em Algebraic Geometry},
\newblock  Springer-Verlag, 1982.

\bibitem{lowa}
K. Stanislaw Lowasiewicz,
\newblock  {\em Introduction to complex Analytic Geometry},
\newblock  Birkhauser Verlag, 1991.

\bibitem{matsumura}
Hideyuki Matsumura,
\newblock  {\em Commutative Ring Theory},
\newblock  cambridge University Press, 1980.

\bibitem{mumford}
David Mumford,
\newblock  {\em The Red Book of Varieties and Schemes},
\newblock  Springer Verlag, 1988.


\bibitem{peternell}
Th. Peternell,
\newblock  {\em Modifications},
\newblock  in Several Complex Variables VII, H. Grauert, Th. Peternell, R. Remmert (Eds), Encyclopedia of Mathematical Sciences {\bf 74} (1994), p. 285-317.

\bibitem{remmert}
R. Remmert,
\newblock  {\em Local Theory of Complex Spaces},
\newblock  in Several Complex Variables VII, H. Grauert, Th. Peternell, R. Remmert (Eds), Encyclopedia of Mathematical Sciences {\bf 74} (1994), p. 7-96.

\bibitem{serre}
J.-P. Serre,
\newblock  {\em Géométrie Algébrique et Géométrie Analytique},
\newblock  Ann. Inst. Fourier, T. {\bf VI} (1956), p. 1-42.

\end{thebibliography}
\end{document}